\documentclass[10pt, a4paper]{article}

\usepackage{fancyhdr}
\usepackage{color}
\usepackage{titletoc}
\usepackage{latexsym}
\usepackage{float}
\usepackage{latexsym,bm}
\usepackage{amsmath}
\usepackage{amssymb}
\usepackage{multicol}
\usepackage{graphics}
\usepackage{graphicx}
\usepackage{indentfirst}
\usepackage{epsfig}
\usepackage{mathrsfs}
\usepackage{algorithm}
\usepackage{algorithmic}
\usepackage{amsthm}
\usepackage{subcaption}
\usepackage{tabularx}

\newtheorem{theorem}{Theorem}

\newtheorem{proposition}{Proposition}

\newtheorem{problem}{Problem}

\def\eop{{\hfill\vbox{\hrule height .3pt
      \hbox{\vrule width.3pt height 7pt
      \kern 7pt
      \vrule width .3pt}
      \hrule height .3pt}} \par\bigskip}

\begin{document}
\title{{\large  \bf Inverse model for network construction: $(\delta(G),I'(G))\to G$}}
\author{{\large  Wei Gao$^{1}$, Yaojun Chen$^{2}$, Hainan Zhang$^{3}$\footnote{Corresponding author: zhnzhanghainan@126.com (H. Zhang).}}
\\
{\small 1. School of Mathematics, Hohai University, Nanjing 210098, China}\\
{\small 2. Department of Mathematics, Nanjing University,
Nanjing 210093, China}\\
{\small 3. School of Software Engineering, South China University of Technology, Guangzhou 510000, China}
\date{}
}

\maketitle
\begin{abstract}
The isolated toughness variant is a salient parameter for measuring the vulnerability of networks, which is inherently related to fractional factors (used to characterize the feasibility of data transmission). The combination of minimum degree and the corresponding tight bound of isolated toughness variant for fractional factors provide reference standards for network construction. However, previous advances only focused on how to select the optimal parameter criteria from Pareto front, without any suggestion for the construction of specific networks. To overcome this deficiency, this paper proposes an inverse model from $(\delta(G),I'(G))$ to $G$, by means of evolutionary computing approach, we propose a novel inverse model to obtain the optimal solutions for candidate graphs. The main procedure is composed of pseudo-greedy acceleration, cross-mutation and diversity enhancement modules. The practicality of the algorithm is verified by means of pilot experiments. The code in this paper is made public on https://github.com/AizhEngHN/Inverse-model-for-network-construction-G-I-G-G.
\end{abstract}



\baselineskip 18pt

\section{Introduction}
The topological structure of a network determines its performance in a geometric sense, and the network structure scheduling during the network construction phase directly influences the network competence and related maintenance (see Stepanova et al. \cite{Stepanova24}, Vagts \cite{Vagts24} and Zarghami \cite{Zarghami24}). After modeling the network by a graph (sites are identified by vertices, and the channels between sites are characterized by edges between vertices), some properties of the network can be measured in terms of graph parameters (see Anitha and Venkataraman \cite{Anitha24}, Khan et al. \cite{Khan23} and Manickam et al. \cite{Manickam20}). Recent works revealed that the combination of minimum degree and isolated toughness parameters provides a reference standard for network designing, and specifically, the isolated toughness variant provides massively distributed parameter combinations (as the minimum degree increases, the tight bound of the isolated toughness variant decreases).

Multi-objective optimization is very instructive for solving conflicting objectives, thus it aims to obtain a set of solutions instead of a unique solution (Kaya and Maurer \cite{Kaya23} and Koziel and Pietrenko-Dabrowska \cite{Koziel22}). Chen et al. \cite{Chen} applied the knee point calculation strategy in multi-objective optimization to select the optimal combination of minimum degree and isolated toughness variant. However, for network designers, knowing the parameter requirements is far from network construction, which is imperative to determine the specific network topology structure, i.e., it entails knowing the specific form of the adjacency matrix of the corresponding graph $G$.

To address this issue, we propose an inverse model that infers graph $G$  (corresponding to the structure of the network) by means of $(\delta(G),I'(G))$ (the specific motivation will be elaborated in the next section, after the comprehensive statement of the related works). To our knowledge, we are the first to consider the network construction problem in terms of inverse model with graph-based parameters. The contributions of this article are delineated below.\\
(1) We propose an inverse model to get the graph structure in light of graph parameter combination $(\delta(G),I'(G))$, which consists of several submodules: initialization, pseudo-greedy acceleration, and diversity enhancement. Furthermore, the time complexity of the given inverse model algorithm is analyzed.\\
(2) To ascertain the availability of the given algorithm, some pilot experiments are presented, and the results signify a greater emphasis on network construction.

The remaining sections are organized as follows. A pragmatic review is manifested in the next section, as well as the problem formulation. Furthermore, we have elaborations  on the motivation and challenges of our work. Next, the proposed evolutionary computation-based algorithm is manifested in Section  \ref{mainalgorithmsection}. Ultimately, to demonstrate the effectiveness of the proposed inverse model algorithm, the pilot experiments are implemented and the obtained results are analyzed.

\section{Related work, motivation and problem formulation}
In this section, preliminary knowledges of graph theory and multi-objective optimization are presented, and we delineate and analyze the previous work in detail. Hence the motivation of our work is elaborated, and formally characterizes the problem ultimately.
\subsection{Notations and terminologies in graph theory}
Only finite undirected simple graphs (networks) are considered in this paper, where $V(G)$ and $E(G)$ are denoted by the vertex set and edge set of a graph $G$ respectively. Let $\delta(G)=\min_{v\in V(G)}\{d(v)\}$ be the minimum degree of $G$, where $d: V(G)\to \mathbb N\cup\{0\}$ is the degree function on $V(G)$. For vertex subset $S$, denote $i(G-S)$ by the number
of isolated vertices in $G-S$ (vertices with degree zero in $G-S$). For two graphs $G_{1}$ and $G_{2}$, we denote $G_{1}\vee G_{2}$ by connecting each pair to vertices from these two graphs. Other standard graph theory notations used but not clearly stated in this paper can follow from \cite{Bondy08}.

Let $a,b,k\in\Bbb N$ with $1\le a\le b$, and $h$ be a fractional indicator function which assigns a real number in [0,1] to each edge.  A {\em
fractional $[a,b]$-factor} is a spanning subgraph determined by the support of $h$ (i.e., $E_{h}=\{e\in E(G)\,|\,h(e)>0\}$) such that
$a\le\sum_{x'\in N(x)}h(xx')\le b$ for any $x\in V(G)$.
We say a graph $G$ admits a fractional $[a,b]$-factor if such $h$ exists.
If $a=b=k$, then it degenerates to fractional $k$-factor. Clearly, if $h$ takes a binary function (only takes value 0 or 1), then the fractional $k$-factor degenerates to the traditional $k$-factor. In computer networks, due to the constraint of channel transmission capacity, large packets need to be partitioned into several small packets, and transmitted through different channels. This problem can be converted to the existence problem of fractional factors, where $[a,b]$ or $k$ characterize the data transmission capacity of the site per unit time.

Yang et al. \cite{Yang01} and Zhang and Liu \cite{zhang01} introduced {\itshape{isolated toughness}} and {\itshape{isolated toughness variant}} which are stated by (if $G$ is not complete)
$$I(G)=\min_{S\subset V(G),i(G-S)\ge2}\left\{\frac{|S|}{i(G-S)}\right\}$$
and
$$I'(G)=\min_{S\subset V(G),i(G-S)\ge2}\left\{\frac{|S|}{i(G-S)-1}\right\},$$
respectively. For a complete graph $G$, it mandates $I(G)=I'(G)=+\infty$ since no vertex subset $S$ meets the restriction $i(G-S)\ge2$.
Given the order of graph $n$ ($n\ge3$), the extreme graph for the minimum value of isolated toughness (variant) is $K_{1,n-1}$ (star network), and the corresponding extreme value is $\frac{1}{n-1}$ (resp. $\frac{1}{n-2}$). As the graph density increases, the values of $I(G)$ and $I'(G)$ increase and tend towards infinity. Meanwhile, it's appreciable noted that for any rational number $\frac{c}{d}$ between $\frac{1}{n-1}$ (resp. $\frac{1}{n-2}$) and $\infty$, a corresponding graph structure can be found with $I(G)=\frac{c}{d}$ (resp. $I'(G)=\frac{c}{d}$) (for instance, $I(K_{c}\vee(dK_{1}))=\frac{c}{d}$ and $I'(K_{c}\vee((d+1)K_{1}))=\frac{c}{d}$). These two parameters are prominent to characterize the stability and vulnerability of the network, the larger the value, the more robust the network, and vice versa.

\subsection{Multi-objective optimization problem (MOP)}
The standard MOP can be formalized by
$$\min {\bf F(x)}=(f_{1}({\bf x}),\cdots,f_{r}({\bf x}))^{T}$$
$${\rm s.t.}\, {\bf x}\in \Omega$$
where ${\bf x}\in\mathbb R^{d}$ is a decision vector, $\Omega$ is a decision space (a set of candidate solutions of MOP), $f_{i}$ ($1\le i\le r$ and $r\ge2$) are objective functions corresponding to $r$ conflicting objectives (in particular, bi-objective optimization problem (BOP) refers to $r=2$). For two potential solutions ${\bf x}^{1},{\bf x}^{2}\in\Omega$, we say ${\bf x}^{1}$ dominates ${\bf x}^{2}$ (denoted by ${\bf x}^{1}\preceq{\bf x}^{2}$) iff $f_{i}({\bf x}^{1})\le f_{i}({\bf x}^{2})$ hold for  $\forall i\in\{1,\cdots,r\}$ and $f_{i}({\bf x}^{1})< f_{i}({\bf x}^{2})$ for at least one $i\in\{1,\cdots,r\}$. A solution ${\bf x}\in \Omega$ is Pareto optimal if there is no other solution ${\bf x}'\in \Omega$ such that ${\bf x}'$ dominates ${\bf x}$, the set of all Pareto optimal solutions is called the Pareto set (PS) and $\{{\bf F(x)}\,|\,{\bf x}\in {\rm PS}\}$ is called the Pareto front (PF).

\begin{figure}[htbp]
    \centering
    \includegraphics[width=1.2\linewidth]{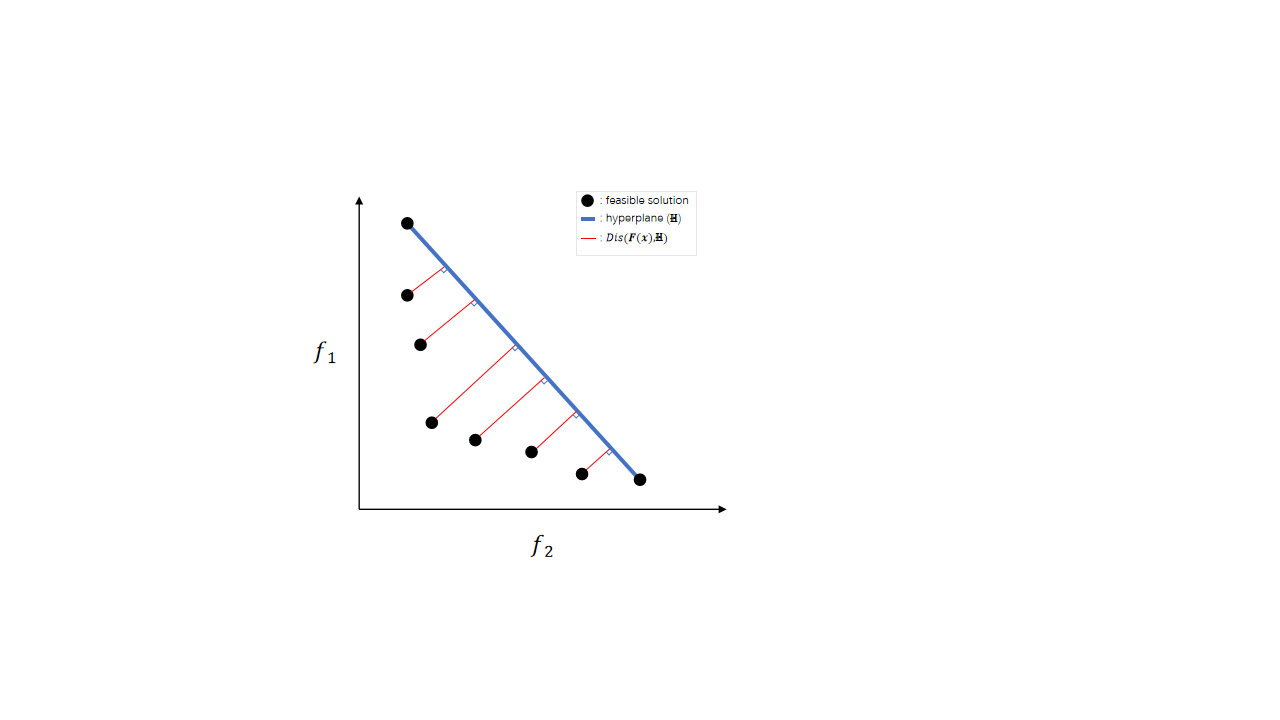}
    \caption{Knee point in BOP}
	\label{KP}
    \end{figure}

The knee point (KP) is the point on the PF with maximum curvature, which can be formulated by
\begin{equation}\label{equ2}
{\bf{x}}_{kp}=\mathop{\arg\max}\limits_{{\bf{x}}\in\Omega}\emph{Dis}({\bf{F(x)},\mathbb{H}}),
\end{equation}
where $\mathbb{H}$ is a hyperplane determined by extreme points on the objective space, and $\emph{Dis}(\cdot)$ is the geometric measure (for example, Euclidean distance). The computing of global knee point in BOP is instantiated in Fig. \ref{KP}.

\subsection{Related work}
From the previous analysis, the larger the value of isolated toughness (variant), the denser the corresponding network and hence the more defensive, but the construction of high-density networks is huge cost and impractical. Network designers are fueled to excavate an equilibrium point where the corresponding value of isolated toughness (variant) cannot be too large or too small. To our delight, theoretical results have revealed an inherent connection between isolated toughness (variant) and fractional factors. The tight isolated toughness (variant) bound of the existence of fractional factors provides a balance point for network designing, enabling the constructed network to satisfy the isolated toughness (variant) boundary which ensures data transmission and network robustness, and on the other hand, controlling the sparsity of the network within a certain range, thereby reducing the construction cost. From this perspective, the tight bound of isolated toughness (variant) for the existence of fractional factors can serve as a parameter criterion for network design.

The sharp $I(G)$ or $I'(G)$ bounds (combined with $\delta(G)$ restriction) for fractional factors in various settings can be summarized in Tab. \ref{tab1}.
\begin{table}
\caption{Tight $(\delta(G),I(G))$ and $(\delta(G),I'(G))$ combination bounds for fractional factors.}\label{tab1}
\begin{tabular}{ccc}
\hline
 Reference & Fractional factor setting & Combine conditions    \\
\hline
Ma and Liu \cite{fractionalkfactor} & fractional $k$-factor & $\delta(G)\ge k$, $I(G)\ge k$ \\
He et al. \cite{He2023} & fractional $k$-factor & $\delta(G)\ge k$, $I'(G)>2k-1$  \\
Gao et al. \cite{abfactor} & fractional $[a,b]$-factor &  $\delta(G)\ge a$, $I(G)\ge a-1+\frac{a}{b}$    \\
Gao and Wang \cite{abfactor2} & fractional $[a,b]$-factor & $\delta(G)\ge a$ and
$I'(G)>a-1+\frac{a}{n_{a,b}-1}$ \\
\hline
\end{tabular}
\end{table}

In the last line of Table \ref{tab1}, $n_{a,b}\ge2$ is an integer such that $(n_{a,b}-1)a\le b\le n_{a,b}a-1$.
It is emphasized that the bounds of the isolated toughness (variants) mentioned in Table \ref{tab1} are all sharp, as evidenced by counterexamples in the corresponding references. More recent results on fractional (component) factors and network parameters can be referred to \cite{recent1} and \cite{recent2}.

Obviously, according to the definition of fractional $k$-factor (resp. fractional $[a,b]$-factor), $\delta(G)\ge k$ (resp. $\delta(G)\ge a$) is tight for fractional $k$-factor (resp. fractional $[a,b]$-factor). The role of minimum degree in \cite{fractionalkfactor}, \cite{He2023}, \cite{abfactor} and \cite{abfactor2} is handling some extreme situations (a dense graph satisfying the $I(G)$ or $I'(G)$ condition, but the minimum degree is less than the lower bound due to small order of graph, and such graph does not admit a fractional factor). Imagine an extreme circumstance: a complete graph with small vertices, its isolated toughness (variant) value is $+\infty$, but without a fractional $k$-factor (assume $k$ is large enough). That is to say, in \cite{fractionalkfactor}, \cite{He2023}, \cite{abfactor} and \cite{abfactor2}, we mainly focus on the tight bound of $I(G)$ (resp. $I'(G)$), while considering the bound of $\delta(G)$ as an auxiliary condition for handling extreme cases, and thus their status are unequal.

Observe that minimum degree and isolated toughness (variant) are two parameters that measure the density of graphs from different perspectives, but they have the same trend (the denser the graph, the larger the corresponding value of minimum degree and isolated toughness (variant); on the contrary, the sparser the graph, the smaller their corresponding values).
Fueled by this prominent fact, Chen et al. \cite{Chen} considered the following problem.
\begin{problem}\label{problem1} {\rm(Chen et al. \cite{Chen})} Considering the combination bound of $\delta(G)$ and $I(G)$ {\rm(}or $I'(G)${\rm)}
for the existence of fractional factors. If the lower bound of the minimum degree is increased, will the tight bound of the corresponding isolated toughness {\rm(}variant{\rm)} decrease?
\end{problem}

For the isolated toughness variant $I'(G)$,  the positive answer is yield by Chen et al. \cite{Chen}, which is described by the following theorem for fractional $k$-factor.
\begin{theorem}\label{mainresult1} {\rm(Chen et al. \cite{Chen})} Let $G$ be a graph, $k\ge2$ be an integer and $t\in\Bbb N\cup\{0\}$. If $\delta(G)\ge k+t$ and $I'(G)>k+\frac{k-1}{t+1}$, then $G$ admits a fractional $k$-factor.
\end{theorem}

And it can be generalized to fractional $[a,b]$-factor setting.
\begin{theorem}\label{mainresult2} {\rm(Chen et al. \cite{Chen})} Let $G$ be a graph, $t\in\Bbb N\cup\{0\}$, and $a,b$ be integers satisfying $2\le a\le b$. If $\delta(G)\ge a+t$ and
$I'(G)>a-1+\frac{a+t}{\lceil\frac{b(t+1)+1}{a}\rceil-1}$,
then $G$ admits a fractional $[a,b]$-factor.
\end{theorem}
By the counterexample presented by Chen et al. \cite{Chen}, the isolated toughness variant bound under strengthened minimum degree constraint is tight.

However, the answer to Problem \ref{problem1} for $I(G)$ is negative (both for fractional $k$-factor and fractional $[a,b]$-factor), i.e., the tight
$I(G)$ bound will not change anymore if increasing the minimum degree bound of the graph. To interpret this phenomenon, we first clarify why in Ma and Liu \cite{fractionalkfactor} and Gao et al. \cite{abfactor}, the tight bounds of $I(G)$ are achievable (``$\ge$''), while in He et al. \cite{He2023} and Gao and Wang \cite{abfactor2}, the sharp bounds of $I'(G)$ are unattainable (``$>$'').

\begin{figure}[th]
\centering
\includegraphics[width=12cm,height=8cm]{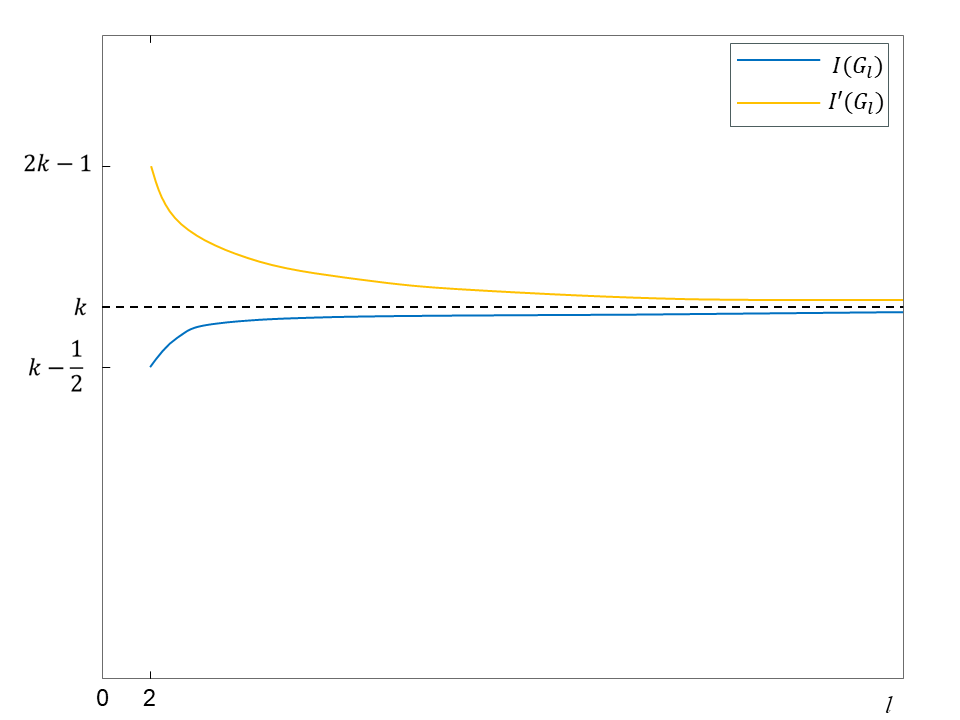}
\caption{Curves for $I(G_{l})$ and $I'(G_{l})$.}\label{figure}
\end{figure}

The underlying mechanisms can be found in Gao et al. \cite{Gao22}, and to be completeness, a detailed explanation is provided below. Consider a graph class $G_{l}=K_{l-1}\vee (lK_{k})$ where $l\in\Bbb N-\{1\}$, and we confirm that $G_{l}$ has no fractional $k$-factor in terms of its sufficient and necessary condition (cf. Liu and Zhang \cite{Liu01}). It easy to check that $\alpha(G_{l})=l$, $I(G_{l})=k-\frac{1}{l}$ and $I'(G_{l})=k+\frac{k-1}{l-1}$ (see Fig. \ref{figure}).

For isolated toughness, $I(G_{l})$ is a monotone increasing function with respect to $l$ and
$\lim_{l\to+\infty}I(G_{l})=k$, and hence explains why tight bound of $I(G)$ in Ma and Liu \cite{fractionalkfactor} is achievable (``$\ge$''). For the isolated toughness variant, $I'(G_{l})$ is a monotone decreasing function with respect to $l$, $\max_{l\ge2}\{I'(G)\}=2k-1$ and it arrives at the maximum value when $l=2$. It explains why the tight bound of $I'(G)$ in He et al. \cite{He2023} is unattainable (``$>$''). In addition,  the intrinsic difference of tight $I'(G)$ bounds between Gao et al. \cite{abfactor} and Gao and Wang \cite{abfactor2} in fractional $[a,b]$-factor setting can be interpreted using the same fashion.

Turning back to original Problem \ref{problem1}, by analysing the proving procedure in Chen et al. \cite{Chen}, we know that increasing the lower bound of $\delta(G)$ will increase the cardinality of the maximum independent set (correspond to $l$ in the foregoing instance). For $I(G_{l})$, its maximum value occurs when $l$ approaches $+\infty$, so finitely increasing the value of $l$  cannot change the lower bound of $I(G)$. On the contrary, for the isolated toughness variant, increasing the cardinality of the maximum independent set leads to an increase in the denominator value of $I'(G)$ of the extreme graph, resulting in decreasing the extreme value of $I'(G)$. This explains why the answer to Question \ref{problem1} regarding the isolated toughness variant is affirmative, but the answer to isolated toughness is negative. It concludes from the aforementioned analyses that $I'(G)$ is dramatically more important than $I(G)$ for network construction, because there is a balance between the $I'(G)$ and the minimum degree, which enables to dynamically turn the graph-based topology parameters of the network from an equilibrium perspective.

In the original theoretical result in He et al. \cite{He2023}, there is only one $(\delta(G),I'(G))$ standard $(k,2k-1)$ for network construction, while after the extension theorem obtained in Chen et al. \cite{Chen}, it becomes infinitely $(\delta(G),I'(G))$ standards for network construction (for $t=0,1,2,\cdots$):\\
$$(k,2k-1),$$
$$(k+1,k+\frac{k-1}{2}),$$
$$(k+2,k+\frac{k-1}{3}),$$
$$\cdots$$
It inevitably leads to a ``choice dilemma'' for network decision-makers, and this parameter selecting problem is dubbed as ``{\em $\delta-I$ equilibrium problem}''. Chen et al. \cite{Chen} regarded decision vector ${\bf x}$ as a graph structure, objective functions $f_ {1}$ and $f_ {2}$ are taken as minimum degree and isolated toughness variant, and $(\delta(G),I'(G))$ is regarded as the PF of a bi-objective optimization problem.
They solved $\delta-I$ equilibrium problem in light of knee point identification (KPI) algorithm, and it is imperative to consider two cases: the vertex number is given or the order of the graph is unknown. Recall the well-known minimum degree condition for fractional $k$-factor which is manifested as follows.
\begin{theorem}\label{mainresult3} {\rm (Gao et al. \cite{Gao20})} Let $G$ be a graph with order $n$, and $k\ge 2$ be an integer. If $n\ge 4k-5$ and $\delta(G)\ge\frac{n}{2}$, then $G$ has a fractional $k$-factor.
\end{theorem}

By the counterexample given by Gao et al. \cite{Gao20}, both bounds for $n$ and $\delta(G)$ are sharp. Hence, when the order of the graph is given, then the hyperplane in objective space is uniquely determined by the extreme points according to the following two circumstances:\\
$\bullet$ If $n\ge4k-5$, then consider $\delta(G)\in[k,\frac{n}{2})$ in Pareto front;\\
$\bullet$ If $n<4k-5$, then consider $\delta(G)\in[k,n-1]$ in Pareto front.\\
For the former case, ``$\frac{n}{2}$'' is called ``transformation point'', which means if the minimum degree excess $\frac{n}{2}$, two conditions $\delta(G)$ and $I'(G)$ is degenerated to a single minimum degree condition, and isolated toughness variant condition no longer functioning (i.e., in this case, the bi-objective optimization problem degenerates into a single-objective optimization problem). However, if the order of the graph is unknown, then the hyperplane in objective space can't be determined since one extreme point trends to infinite when $t\to\infty$. In this case a novel KPI algorithm (see Chen et al. \cite{Chen}) is tremendously used to approximate the knee point which to be selected as the parameter criterion in network designing.

\subsection{Motivation and challenge}
Although \cite{Chen} provides a selection strategy for network construction parameter standards, for network engineers (decision makers), knowing the parameter criterion is far from real engineering requirements. Determining the actual network structure (the adjacency matrix of $G$) is paramount in many application tasks. Suppose $G$ is an undirected graph with $n$ vertices, then the dimensionality of ${\bf x}$ in decision space is ${n \choose 2}$ which determines whether there is an edge between each pair of vertices in $G$. This inspires us to consider obtaining the neighborhood structure of $G$ from two-dimensional data $(\delta(G),I'(G))$, and it is equivalent to solving an inverse model from low dimensional space to high-dimensional space.

To construct the above-mentioned inverse model faces the following challenges:\\
(1) The isolated toughness variant is a nonlinear function that cannot be characterized in matrix form, and therefore the folklore inverse models in MOP are failures in our setting.\\
(2) Given a graph $G$, then $\delta(G)$ and $I'(G)$ are uniquely determined. But conversely, given the values of $\delta(G)$ and $I'(G)$, there are many graphs that satisfy the combination condition. We need a selection strategy for these ${\bf x}$ that meets the conditions.\\
(3) For the construction of the initial population, we naturally think of the counterexample $K_{t+1}\vee((t+2)K_{k})$ in Chen et al. \cite{Chen}. However, the following points should be fully considered when using these counterexamples as individuals.\\
$\bullet$ First, since $\delta(G)\in\Bbb N$, the individuals only distributed on vertical lines correspond to the integer values of the minimum degree axis.\\
$\bullet$ Second, there is no fractional $k$-factor in the counterexample, so theoretically the points located on the PF cannot be taken as individuals. In other words, if the PF corresponding to the counterexample is used as the boundary of individual set, then it is an open set, and the points on the boundary cannot be taken. From a theoretical point of view, the ideal individuals are points on the integer vertical line that is relatively close to the PF boundary.\\
 $\bullet$ Thirdly, the points of graphs corresponding to counterexample $K_{t+1}\vee((t+2)K_{k})$ are located exactly on the PF, but it does not mean that all points on the PF belong to this graph class. There may be other families of graphs that meet the numerical characteristics of the counterexample (i.e., $\delta(G)=k+t$ and $I'(G)=k+\frac{k-1}{t+1}$).\\
$\bullet$ Finally, in network construction, the value $n$ is fixed which embodies the number of sites in the network, and $k$ and $t$ are both integers. It is very likely that we cannot find an integer combination $(n,k,t)$ satisfying $n=t+1+k(t+2)$. This makes it impossible to achieve the goal of constructing individuals through counterexamples, and thus it is necessary to relax the conditions, such as allowing the value of $k$ to be slightly larger than the predefined configuration.

We end this section by showcasing Tab. \ref{tab2} which lists all the appreciable symbols appearing in the paper and provides corresponding descriptions.
\begin{table}[!htp]
    \centering
    \caption{Notations} \label{tab2}
    \begin{tabularx}{\textwidth}{c>{\centering\arraybackslash}X}
        \hline
        Notations  & Description  \\
        \hline
$V(G)$ & The vertex set of graph $G$.\\
$E(G)$ & The edge set of graph $G$.\\
$d(v)$ & The degree of vertex $v$.\\
$\delta(G)$ & The minimum degree of graph $G$  \\
$\Delta(G)$ & The maximum degree of graph $G$.  \\
$I(G)$ & Isolated toughness.  \\
$I'(G)$ & Isolated toughness variant.  \\
$ k$ & The scale of fractional factor.\\
$n$ & The number of vertices in the graph (the order of graph).\\
$t$ & A non-negative integer parameter in the counterexample.\\
$S$ & The vertex subset of the graph.\\
$i(G-S)$ & The number of isolated vertices remaining in graph $G$ after the removal of its subset $S$.\\
\hline
$k_{min}$ & The lower bound of minimum degree. \\
$P$ & The population.\\
$N$ & The size of population $P$.\\
$r$ & The number of objectives.\\
$\Omega$ & The decision space.\\
${\bf x}$ & A decision vector (an individual in the population).\\
${\bf x}_{kp}$ & An individual corresponds to global knee point.\\
$D$ & Individual dimension.\\
$f_{i}$ & Objective function.\\
$\mathbb{H}$ & A hyperplane in objective space.\\
$dis(\cdot,\cdot)$ & Distance function in objective space.\\
$Tmax$ & Maximum number of iterations.\\
${\bf{Deg}}$ & The set of degrees for each vertex.\\
${\textbf{\emph{Ix}}}$ & A representation of an individual {\bf x}.\\
$\mu$ & The probability of binary mutation.\\
$r_{1}$ &  When there is a counterexample population, the proportion of the first part individuals in the total population.\\
$p$ &  A real number in $[0,1]$ which used in roulette wheel strategy.\\
$\mathcal{G}$ or $\mathcal{G'}$ & A set of graphs (a family of networks). \\
$\Theta$ & An intermediate set to store graphs that meet the conditions.\\
$\Theta_{i}$ & A subset of $\Theta$ with minimum degree $i$.\\
\hline
\end{tabularx}
\end{table}

\section{Proposed algorithm}\label{mainalgorithmsection}
The aim of this section is to provide our inverse model algorithm in detail. The main procedure consists of several modules, and it conducts a detailed analysis of each module which is separated in subsequent subsections. It is noteworthy that the algorithm proposed in this article only focuses on graph parameters with respect to fractional $k$-factors (corresponds to Theorem \ref{mainresult1}). Since fractional $[a,b]$-factor can be seen as a relaxation of fractional $k$-factor, the proposed inverse model algorithm can be applied to fractional $[a,b]$-factor setting by virtue of certain relaxations.

\subsection{Model formulation}
\begin{figure}[thb]
\centering
\includegraphics[width=0.8\linewidth]{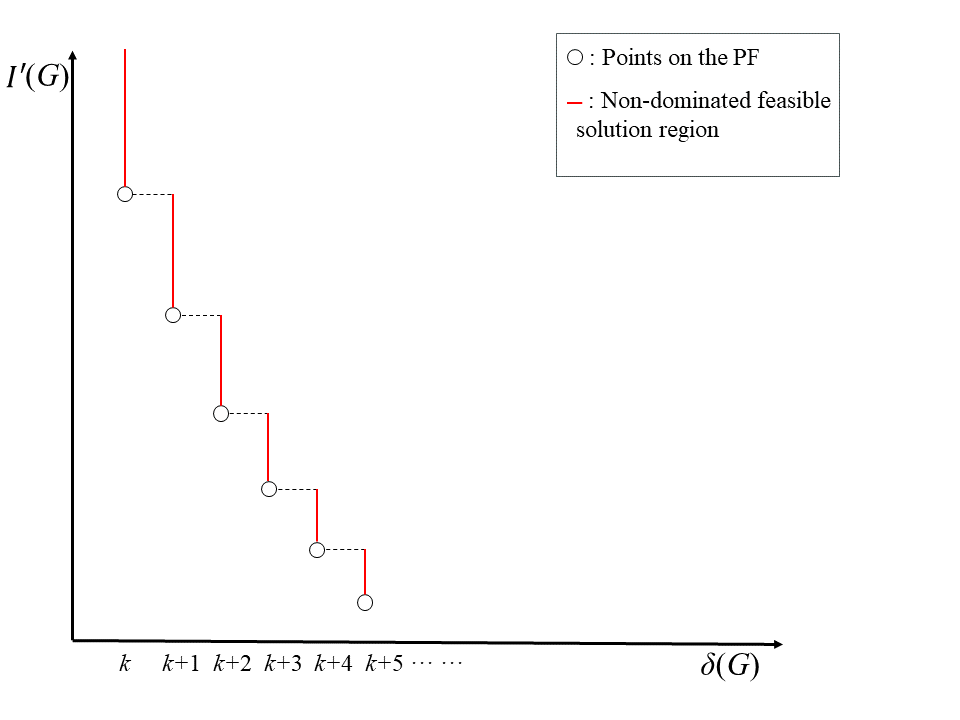}
\caption{Schematic diagram of the objective space.}\label{fig3}
\end{figure}

According to the requirements of network constructors, we initially determine the number of vertices $n$ and the lover bound of minimum degree $k_{min}$ for the network to be constructed, where $n$ is immutable, while $k_{min}$ serves as a lower bound of $\delta(G)$ of a specific network $G$. Hence, a lower bound of theoretically optimal solutions can be derived in the objective space to form a PF. It is important to note that the constructed PF is unattainable since $I'(G)>k+\frac{k-1}{t+1}$.  A schematic diagram of the objective space is  instantiated in Fig. \ref{fig3}, where the two objectives are labeled as $\delta(G)$ and $I'(G)$, respectively. Due to Theorem \ref{mainresult3}, we have an upper bound for $\delta(G)$ when $n$ is given.

The black hollow circles represent points on the PF, which embody they are unattainable. The red solid lines denote feasible non-dominated solutions. Since $\delta(G)$ consists of positive integers, the lines appear as discontinuous solid lines in the objective space. It is emphasized that since the set of rational number is dense in real numbers, without causing confusion, each vertical line in Fig. \ref{fig3} is drawn by a solid line. However, in reality, since $I'(G)$ only take rational numbers, the representation of Fig. \ref{fig3} is not rigorous, just for convenient expression. The specific problem modeling can be formulated by the following equation:
\begin{equation}\label{equMOP}
\begin{split}
&\min {\textbf{\emph{F}}}(\textbf{x}) = [f_1({\bf x}),f_2({\bf x})]\\
&{\rm s.t.}\quad  \left\{\begin{array}{lc}
f_1({\bf x}) = k+t \in \mathbb {N}\\
f_2({\bf x})>k+\frac{k-1}{t+1} \,{\rm and}\, f_2({\bf x})\in\mathbb {Q},\\
\end{array}\right.
\end{split}
\end{equation}
where ${\bf x}=[x_1,x_2,\cdots,x_D]$ is the decision vector which represents a specific graph structure, $D=\frac{n(n-1)}{2}$ represents the length of the decision vector, and $f_1(\textbf{x})$ and $f_2(\textbf{x})$ correspond to $\delta(G)$ and $I'(G)$ respectively.

\subsection{Overall algorithm process}

\begin{figure}[thb]
\centering
\includegraphics[width=1\linewidth]{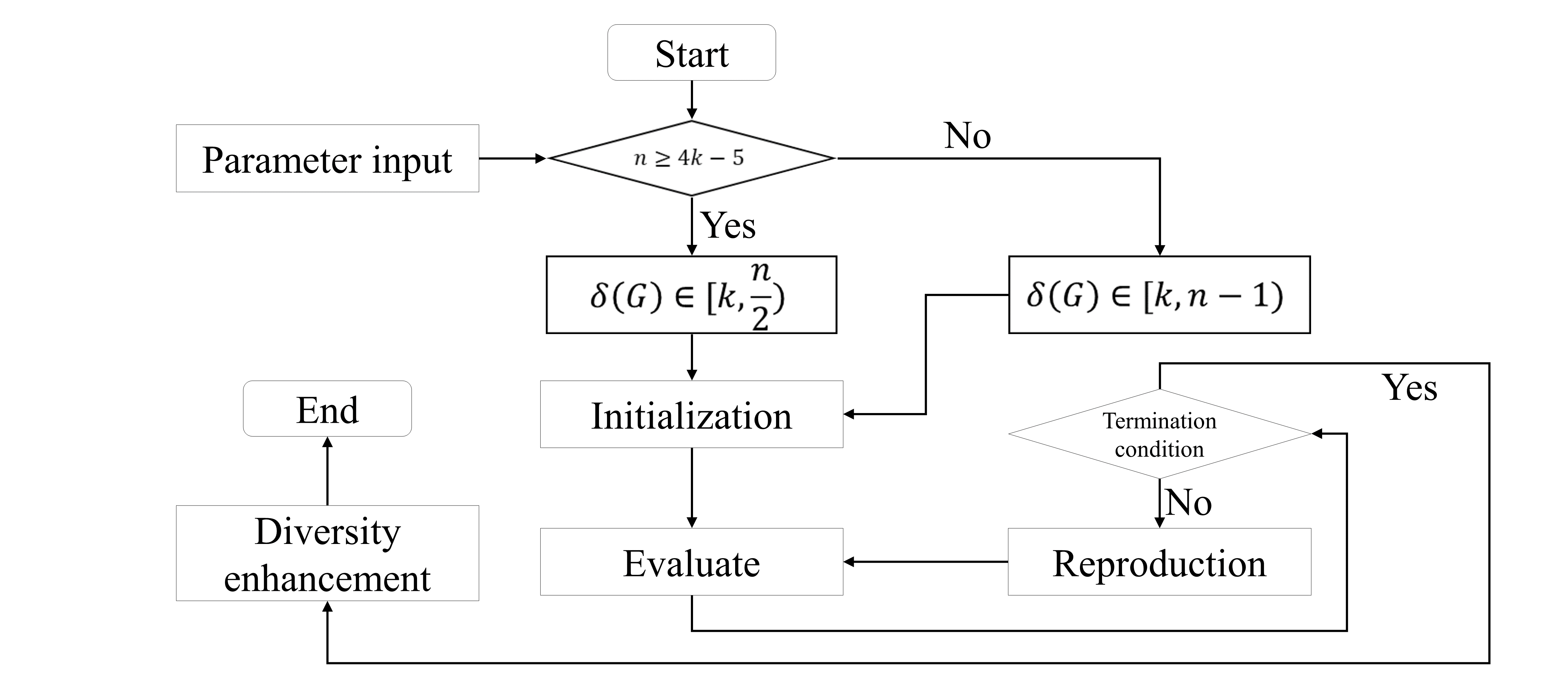}
\caption{The flowchart of main inverse model algorithm.}\label{flew}
\end{figure}

The algorithm we propose consists of several key modules: initialization, pseudo-greedy acceleration, cross-mutation, and diversity enhancement. Alg. \ref{alg1} demonstrates the pseudocode of the algorithm workflow and the overall algorithm implementation is showcased in Fig. \ref{flew}.

\begin{algorithm}
	\caption{Overall algorithm process $P$}
	\label{alg1}
	\begin{algorithmic}[1]
		\REQUIRE $n$: Number of vertices, $\delta(G)$: the minimum degree.
		\ENSURE $\mathcal{G'}$: A representative set of graph structures that meet the requirements.
		\STATE Determine the scope of $\delta(G)$.
		\STATE $P$ = {\bf{Initialization}}($N$,$r_1$)
        \STATE $\mathcal{G}=\emptyset$
        \STATE $j=0$
		\STATE {\bf {while}} Not reaching the termination condition ($j<Tmax$).
        \STATE \quad\quad $\Theta=\emptyset$
        \STATE \quad\quad{\bf{for}} ${\bf x}\in P$
		\STATE\quad\quad\quad\quad  {\bf{Pseudo-greedy acceleration(x)}}.
         \STATE\quad\quad\quad\quad $\Theta=\Theta\cup\{{\bf x}\}$ if ${\bf x}$ meets the requirements.
		\STATE \quad\quad{\bf{end for}}
         \STATE \quad\quad $\Theta=\Theta_{\delta(G)}\cup\Theta_{\delta(G)+1}\cup\cdots\cup\Theta_{\Delta(G)}$ where $i\in\{\delta(G),\cdots,\Delta(G)\}$ and $\Theta_{i}$ is the subset of $\Theta$ whose element has minimum degree exactly equals to $i$.
		\STATE \quad\quad{\bf{for}} $i$ from $\delta(G)$ to $\Delta(G)$
		\STATE\quad\quad\quad\quad  ${\bf x}^{ji}$ $\gets$ Select the individual from  $\Theta_{i}$ with smallest $I'(G)$.
        \STATE\quad\quad\quad\quad  $\mathcal{G}\gets \mathcal{G}\cup\{{\bf x}^{ji}\}$
        \STATE\quad\quad\quad\quad  $i=i+1$
		\STATE \quad\quad{\bf{end for}}
		\STATE\quad\quad $P$ = {\bf{Cross-mutation($P$)}}
        \STATE\quad\quad $j=j+1$
		\STATE {\bf {end while}}
		\STATE $\mathcal{G'}$={\bf {Diversity enhancement}}($\mathcal{G}$)
		\STATE return $\mathcal{G'}$
	\end{algorithmic}
\end{algorithm}

\subsection{Initialization}

First, we initialize a population $P$ to prepare for the evolutionary process, where $P$ is divided into two parts. The first part consists of individuals generated based on the coordinates of the hollow points in Fig. \ref{fig3}. These individuals do not meet the requirements, therefore we apply binary mutation \cite{BM} to fine-tune their encoding. The second part is composed of randomly generated individuals, which aim to fill the population quota and enhance the diversity of the population. Alg. \ref{alg2} describes the procedure of population initialization.

\begin{algorithm}

	\caption{Initialization $P$}
	\label{alg2}
	\begin{algorithmic}[1]
		\REQUIRE $n$: The order of graph; $N$: Population size; $r_1$: The proportion of the first part, where $r_{1}\in(0,1)$.
		\ENSURE $P$
		\STATE $i=0$
		\STATE {\bf {while}} $i \ \textless \ N \times r_1$
		\STATE\quad\quad  $P_1 \gets K_{t+1} \vee ((t+2)K_k)$  // The first part of the population is generated based on counterexamples.
		
		\STATE\quad\quad $i = i + 1$
		\STATE {\bf {end while}}
		\STATE $P_1$ = {\bf{Binary mutation($P_1$)}}
		\STATE $P_2 \gets$ Randomly generate $(1-r_1)N$ individuals.
		\STATE  $P$ = $P_1 \cup P_2$

	\end{algorithmic}
\end{algorithm}

Since the parameters of the counterexample are all integers, the counterexample only exists under certain conditions, i.e., $n=t+1+(t+2)k$. In practical applications, $n$ represents the number of sites, which is fixed in the early stage of network scheduling. $\delta(G)$ describes the lower bound of the channel connecting ability for each site. For a specific individual, the value of the minimum degree can break through the given lower bound. The parameter $k$ characterizes the lower bound of the data transmission capacity of each site per unit time, and hence the data transmission capacity larger than $k$ is allowed for a specific individual. Furthermore, $I'(G)$ elaborates the lower bound of the network's defense capability and is also varying for a specific individual. The physical meaning of these parameters in practical applications can be found in Table \ref{physicalmeaning}.
\begin{table}[!htp]
    \centering
    \caption{Physical meaning of notations in network setting} \label{physicalmeaning}
    \begin{tabularx}{\textwidth}{>{\centering\arraybackslash}X>{\centering\arraybackslash}X>{\centering\arraybackslash}X}
        \hline
        Notation  & Physical meaning & Is the actual value fixed for individuals? \\
        \hline
        $n$ & The number of sites. & $\surd$ \\
        $\delta(G)$ & The lower bound of the channel connecting ability for each site. & $\times$ \\
        $k$ & The lower bound of the data transmission capacity of network. & $\times$ \\
	   $I'(G)$ & The lower bound of the network's defense capability. & $\times$\\
        \hline
    \end{tabularx}
 \end{table}
It uncovers that for each individual, the values of $\delta(G)$, $k$ and $I'(G)$ can be larger than the pre-set lower bound, which allows us to generate new individuals that meet the requirements by transforming counterexamples in terms of crossover and mutation.

\begin{figure}[thb]
\centering
\includegraphics[width=0.8\linewidth]{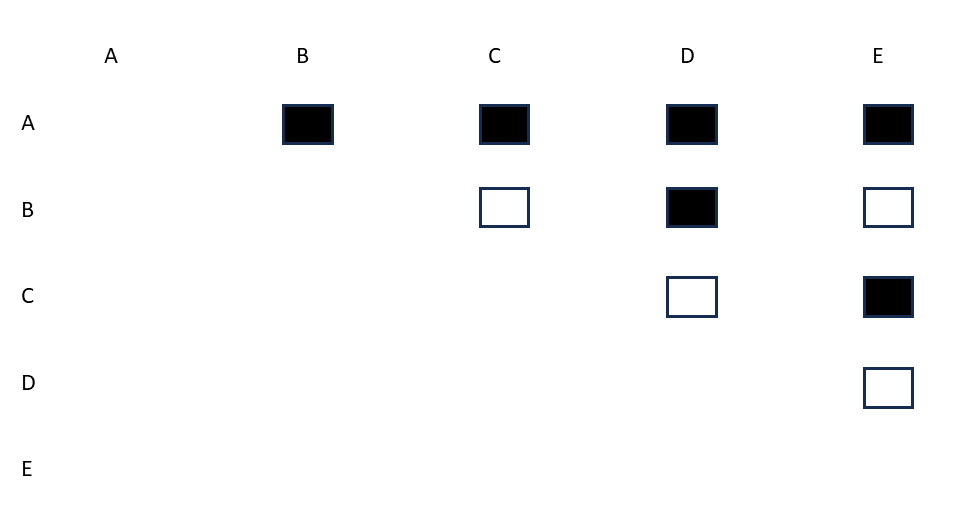}
\caption{Representation of an individual.}\label{fig4}
\end{figure}

Each individual ${\textbf{\emph{Ix}}}$ is represented by a pair of elements, ${\textbf{\emph{Ix}}} = [{\bf{Deg}},{\textbf{\emph{x}}}]$, where ${\bf{Deg}}$ is a vector of degrees for each vertex, and ${\textbf{\emph{x}}}$ is the decision vector, which is a graph structure. For example, assume $n=5$, Fig. \ref{fig4} is a schematic diagram of an individual, where solid dots represent edges between two vertices, and hollow dots indicate no edges. The representation of this individual is denoted by ${\textbf{\emph{Ix}}} = [[4,2,2,2,2],[1,1,1,1,0,1,0,0,1,0]]$.

\subsection{Pseudo-greedy acceleration}
The most toughest step in constructing a network structure as desired by the decision-maker is the computation of $I'(G)$, which is a process with very high computational complexity.  Therefore, conventional execution ideas for network construction are extremely time-consuming. In this paper, we propose a pseudo-greedy acceleration method to expedite the computation of $I'(G)$. We first consider the vertex with the maximum degree in the graph structure, where a larger degree implies more edges incident to this vertex. Removing this vertex can maximally disrupt the graph structure, resulting in a sparser graph with a large number of (potential) isolated vertices. However, it is worth noting that deleting the vertex with the maximum degree does not necessarily result in the most isolated vertices.

\begin{figure}[thb]
\centering
\includegraphics[width=0.8\linewidth]{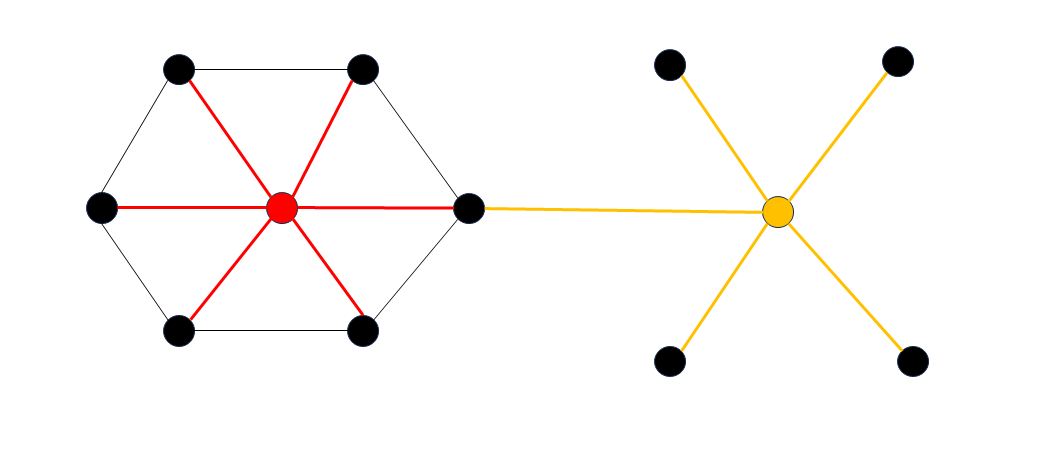}
\caption{Removing the vertex with the maximum degree does not necessarily yield the most isolated points.}\label{fig5}
\end{figure}

As shown in Fig. \ref{fig5}, the red vertex represents the vertex with the maximum degree in the graph. However, removing it does not result in any isolated vertices, as all other vertices are associated with additional edges. On the other hand, removing the yellow vertex in the graph leads to the creation of four isolated vertices. This implies that we cannot solely focus on the vertex with the maximum degree. Therefore, in this paper, in addition to the vertex with the highest degree of deletion, we also adopted a pseudo-greedy acceleration based on the roulette wheel method. The higher the degree of the vertex, the greater the probability of being selected. The specific formula is as follows:
\begin{equation}\label{lunpan}
\frac{\sum^{j}_{i=1} {Deg_i}}{\sum^{\left| {\bf{Deg}} \right|}_{i=1} {Deg_i}} \leq p \textless \frac{\sum^{ j+1 }_{i=1} {Deg_i}}{\sum^{\left| {\bf{Deg}} \right|}_{i=1} {Deg_i}},
\end{equation}
where $p \in [0,1]$ is a random number, $j$ represents the position of the vertex to be removed. The pseudocode for ${\bf{Pseudo-greedy acceleration}}$ is shown in Alg. \ref{alg:pga}.

\begin{algorithm}

	\caption{Pseudo-greedy acceleration}
	\label{alg:pga}
	\begin{algorithmic}[1]
		\REQUIRE ${\bf{x}}$: The individual
		\ENSURE $I'(G)$
		\STATE $s$ = 0 // The number of deleted vertices.
\STATE $I'(G)s$ = 9999 // Initialize to a large number.
\STATE ${\bf{x_I}}$ =  ${\bf{x}}$ , ${\bf{x_{II}}}$ =  ${\bf{x}}$
		\STATE {\bf {for}} $l$ from 1 to $n-3$
\STATE\quad\quad $s$ = $s+1$
		\STATE\quad\quad  $I \gets $  Randomly generate a number from 0 to 1 and use Equ. \ref{lunpan} to select vertex for ${\bf{x_I}}$. // Vertex to be deleted is determined using the roulette wheel selection method.
		\STATE\quad\quad  $II \gets argmax{\bf{Deg(x_{II})}}$  // Vertex to be deleted is determined which has the max degree.
\STATE\quad\quad  ${\bf{x_I}}$ =  ${\bf{x_I}} -I$ , ${\bf{x_{II}}}$ =  ${\bf{x_{II}}} -II$
		\STATE\quad\quad ${\bf {if}}$ $i({\bf{x_I}} ) \geq 2$
\STATE\quad\quad\quad\quad ${\bf {if}}$  $I'(G) \textgreater I'({\bf{x_I}})$
\STATE\quad\quad\quad\quad \quad\quad $I'(G) = I'({\bf{x_I}})$
\STATE\quad\quad\quad\quad ${\bf {else}}$
\STATE\quad\quad\quad\quad\quad\quad ${\bf{x_I}}$ =${\bf{x_{II}}}$
\STATE\quad\quad ${\bf {else}}$
\STATE\quad\quad\quad\quad ${\bf{x_I}}$ =${\bf{x_{II}}}$
		\STATE\quad\quad ${\bf {if}}$ $i({\bf{x_{II}}}) \geq 2$
\STATE\quad\quad\quad\quad ${\bf {if}}$  $I'(G) \textgreater I'({\bf{x_{II}}})$
\STATE\quad\quad\quad\quad \quad\quad $I'(G) = I'({\bf{x_{II}}})$
		\STATE {\bf{end for}}
\STATE return $I'(G)$

	\end{algorithmic}
\end{algorithm}

\subsection{Diversity enhancement}

Due to the mapping from decision space to objective space not being one-to-one but rather many-to-one, there could be many graph structures that satisfy the same requirements. Therefore, it is important to provide decision-makers with as diverse graph structures as possible, rather than offering similar ones. To meet this requirement, we employ a method called diversity enhancement, which aims at increasing the diversity of solutions. The core idea of this trick is to maximize the distance between solutions.

Firstly, since every individual in the final result is a feasible solution, and its encoding represents the upper triangular part of a graph adjacency matrix, we entail eliminating graph structures with the same meaning (isomorphism graphs or graphs with high similarity). Subsequently, we select a graph that has the maximum Hamming distance from the complete graph as the first chosen graph structure. Finally, from the remaining graphs, we select those that have the maximum Hamming distance from the already chosen graph structures. This can be expressed using the following formula:
\begin{equation}\label{DE}
G=\mathop{\arg\max}_{G\in \mathcal{G}-\mathcal{G'}} Dis(\mathcal{G'},G),
\end{equation}
where $G$ represents the chosen graph, $\mathcal{G}$ represents the feasible solutions found by the algorithm, $\mathcal{G'}$ is a representative set of graph structures that meet the requirements, and $Dis(\cdot,\cdot)$ represents the function of Hamming distance. The specific implementation process of diversity enhancement is shown in Alg. \ref{alg3}.

\begin{algorithm}

	\caption{Diversity enhancement}
	\label{alg3}
	\begin{algorithmic}[1]
		\REQUIRE $\mathcal{G}$: Feasible solutions found by the algorithm.
		\ENSURE $\mathcal{G'}$:  A representative set of graph structures that meet the requirements.
		\STATE Eliminate individuals with duplicate graph structures (isomorphism graphs).
\STATE $\mathcal{G'}=\emptyset$
		\STATE $\mathcal{G'} \gets$ Find the graph with the maximum Hamming distance from $K_n$ in $\mathcal{G}$.
\STATE {\bf {while}} $|\mathcal{G'}| \textless N$
		\STATE\quad  $G\gets$ Find the graph in $\mathcal{G}-\mathcal{G'}$ that has the maximum Hamming distance from $\mathcal{G'}$. // Execute (\ref{DE})
       \STATE\quad $\mathcal{G'}\gets \mathcal{G'}\cup\{G\}$
\STATE {\bf {end while}}
		\STATE return $\mathcal{G'}$

	\end{algorithmic}
\end{algorithm}

\subsection{Time complexity analysis}
The overall time complexity of the algorithm is composed of the following parts: initialization, evolutionary process, and diversity enhancement. The time complexity of the initialization part is divided into two scenarios: the first scenario is completely random, which with time complexity $O(Nn)$. The second scenario is partially random and partially utilizes counterexamples for generation, which costs $O(rNn + (1-r)Nn)$. Within the evolutionary process, the time complexity of a single iteration is divided into evaluation, reproduction, and selection, which are performed linearly, resulting in $O(Nn + N^2 + kN)$. The time complexity for diversity enhancement is $O(N^2)$. Therefore, the total time complexity of Alg. \ref{alg1} is $O(Tmax(Nn+N^2+kN))$.

\section{Pilot experiment}
Since the current research does not involve constructing specific graph structures by means of graph-based parameters (to our knowledge, this paper is the first work to construct a network in terms of graph-based parameters, specifically using $\delta(G)$ and $I'(G)$),  we do not compare our algorithm with other methods. Instead, we demonstrate the feasibility of our algorithm with the results from some pilot experiments.

\subsection{Parameter settings}
Since all individuals are undirected graphs, we use binary encoding to represent an individual, with its dimension $D=\frac{n(n-1)}{2}$. During the population evolution process, we use single-point crossover and binary mutation to update the population for binary-encoded individuals. Below, we provide the specific settings for all relevant parameters in this experiment.\\
$\bullet$ The population size $N=10$.\\
$\bullet$ The maximum number of iterations $Tmax=100$.\\
$\bullet$ The probability of binary mutation $\mu=0.3$.\\
$\bullet$ When there is a counterexample population, the proportion of the first part of population is set to $r_{1}=0.5$.

It's noteworthy that since counterexamples are constructed based on $K_{t+1}\vee((t+2)K_k)$ where $t$ must be a positive integer, it is possible that a graph corresponding to the counterexample may not exist when both $n$ and $k$ are fixed. In this case, all individuals are randomly generated.

All experiments in this paper were conducted on a Windows 10 system using Python 3.9. The hardware specifications include a 5700X 8-core CPU and 8GB of RAM.

\subsection{Experimental result}
Subsequently, we shall delineate several sets of encodings corresponding to the optimal graph structures, as derived by the algorithm for specified parameters $n$ and $k$. These encodings serve as a representation of the algorithm's efficacy in determining graph configurations under given constraints. Tab. \ref{tab3} presents the results of the minimum degree and isolated toughness variant obtained for a given number of vertices $n$ and minimum degree $k$, while Tab. \ref{tab4} showcases the number of graph structures obtained (average and variance of thirty experiments) before the implementation of diversity enhancement in the algorithm.

\begin{table}[!htp]
    \centering
    \caption{The results of ($\delta(G)$, $I'(G)$) given different numbers of vertices $n$ and initial minimum degree $k$.} \label{tab3}
    \begin{tabular}{ccllllll}
        \hline
        $n$  & $k$& \multicolumn{6}{c}{ $(\delta(G),I'(G))$}  \\
        \hline
			7		 &      2         &         (2,5)  & (3,5)    &  &  &  &     \\
			9		 &      2         &         (2,7)  &  (3,3)  & (4,3)  &  &  & \\
			11		 &      2         &         (2,4)  &  (3,4) &  (4,4) &  (5,2.333)    &  &       \\
			13		 &      2         &        (2,5) &  (3,3) &  (4,3) &  (5,3) &  (6,3) &    \\
			15		 &      2         &        (2,3.666) &  (3,3.666) &  (4,2.5) &   (5,2.5) &  (6,2.5) & (7,3.666)   \\
			15		 &      3         &         (3,6.0) & (4,6.0) &  (5,6.0) &  (6,3.666) &  (7,3.666) &  \\
			15		 &      4         &         (4,Null) &  (5,6.0) &  (6,6.0) &  (7,6.0) &  &  \\
        \hline
    \end{tabular}
 \end{table}

\begin{figure}[htbp]
	\centering
	\begin{subfigure}{0.325\linewidth}
		\centering
		\includegraphics[width=1\linewidth]{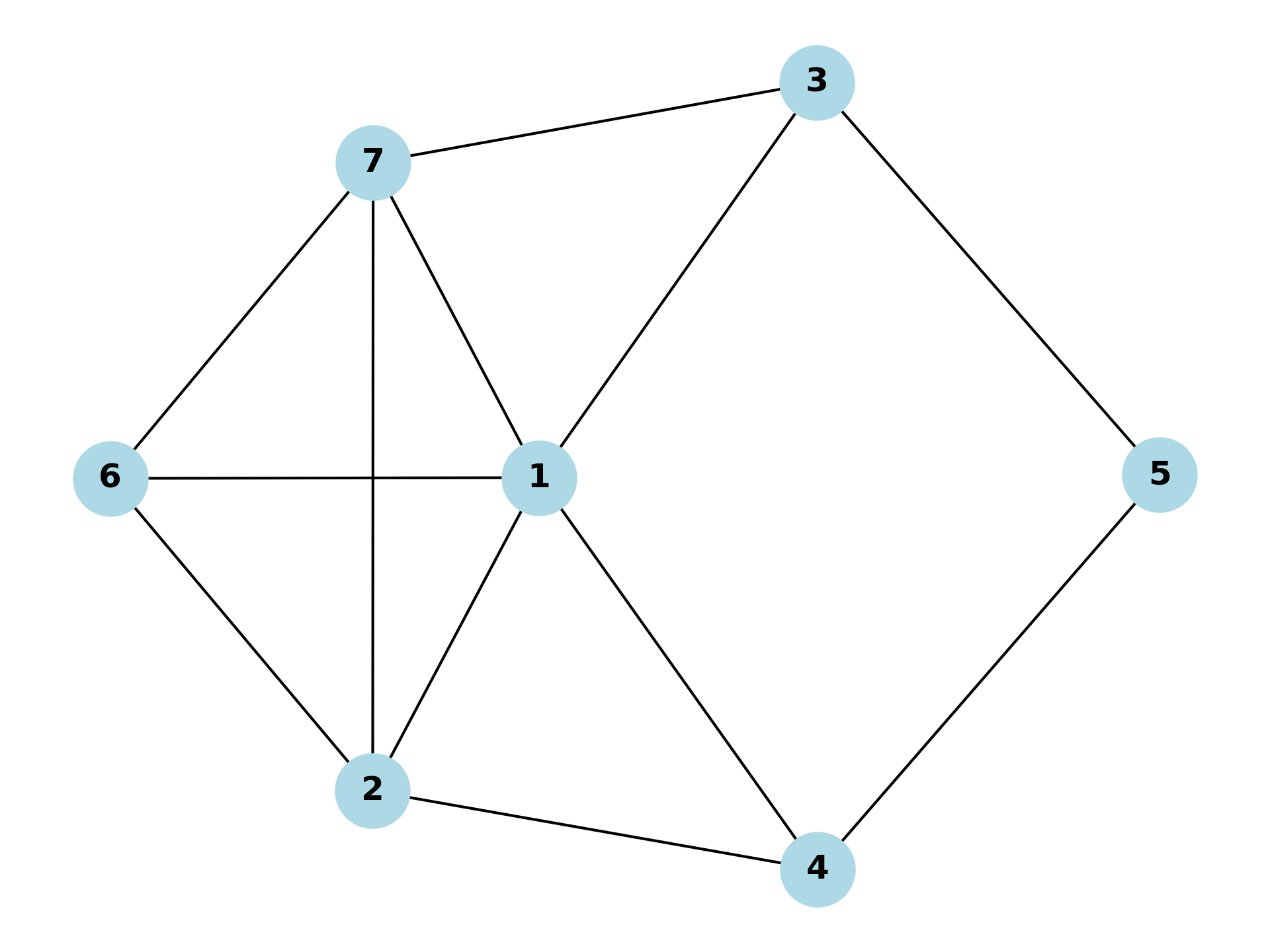}
		\caption{$n=7,\delta(G)=2,I'(G)=5$}
		\label{fig71}%文中引用该图片代��?
	\end{subfigure}
	\centering
	\begin{subfigure}{0.325\linewidth}
		\centering
		\includegraphics[width=1\linewidth]{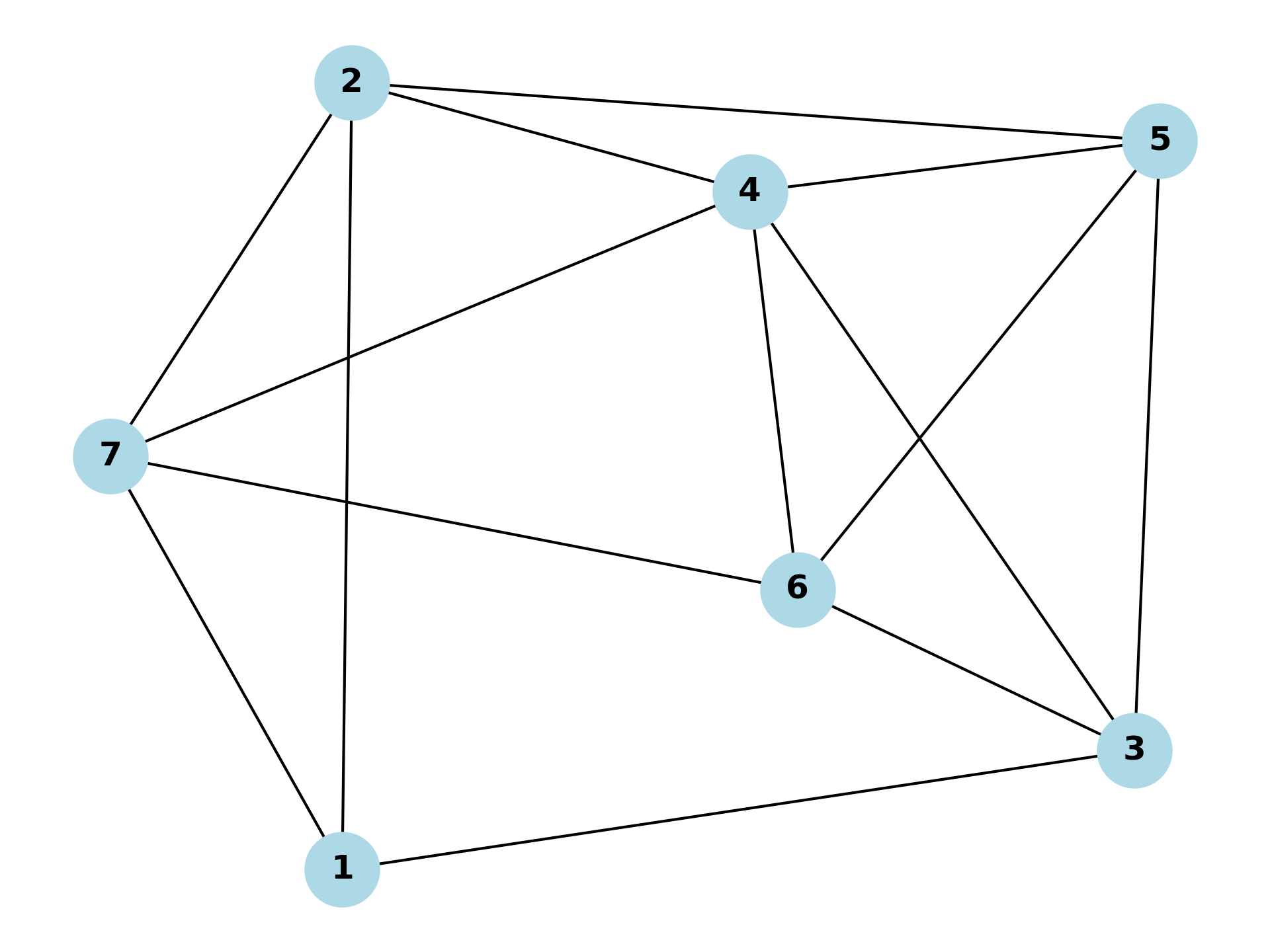}
		\caption{$n=7,\delta(G)=3,I'(G)=5$}
		\label{fig72}%文中引用该图片代��?
	\end{subfigure}
	\centering
	\begin{subfigure}{0.325\linewidth}
		\centering
		\includegraphics[width=1\linewidth]{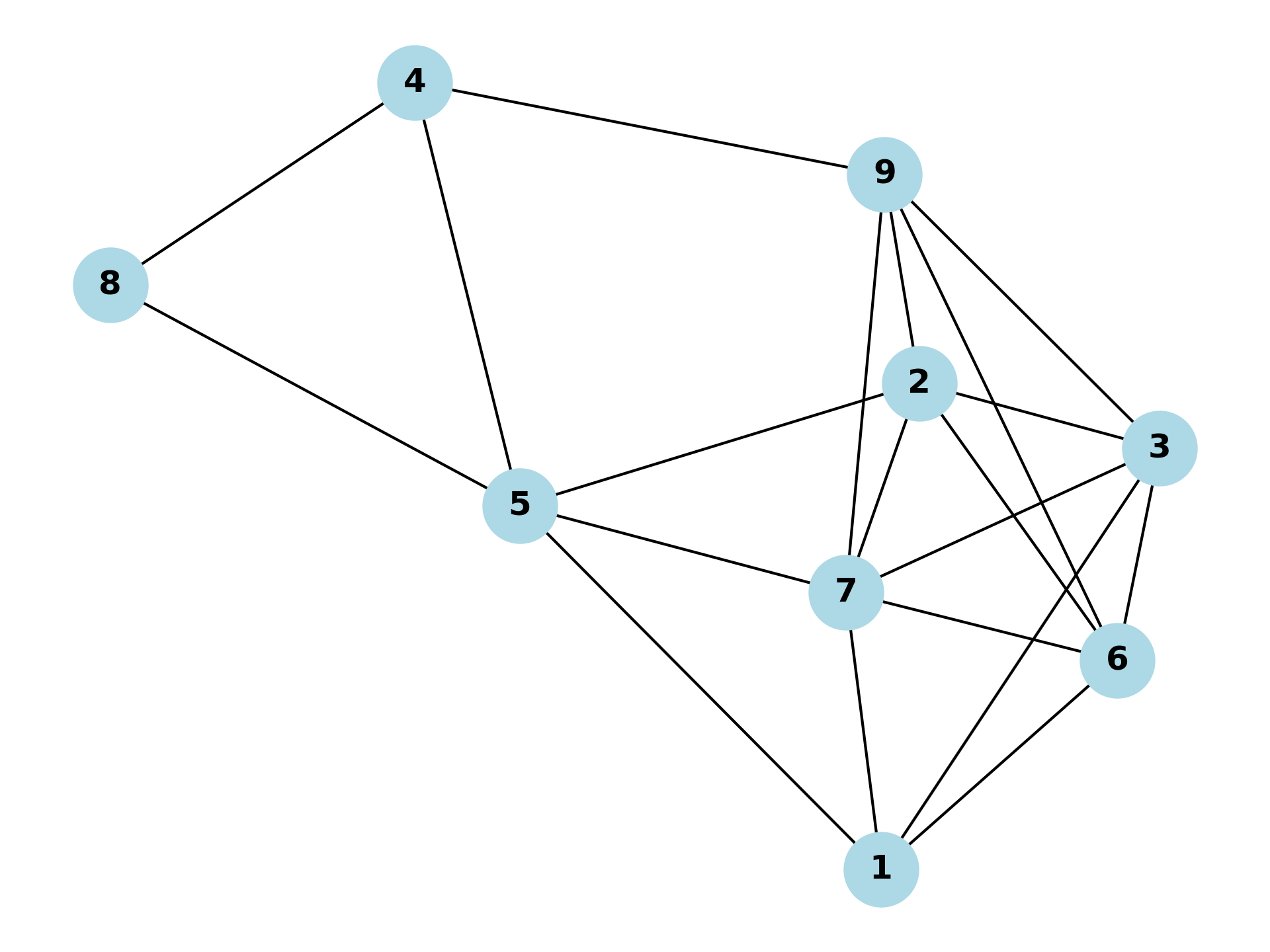}
		\caption{$n=9,\delta(G)=3,I'(G)=7$}
		\label{fig73}%文中引用该图片代��?
	\end{subfigure}
	\begin{subfigure}{0.325\linewidth}
		\centering
		\includegraphics[width=1\linewidth]{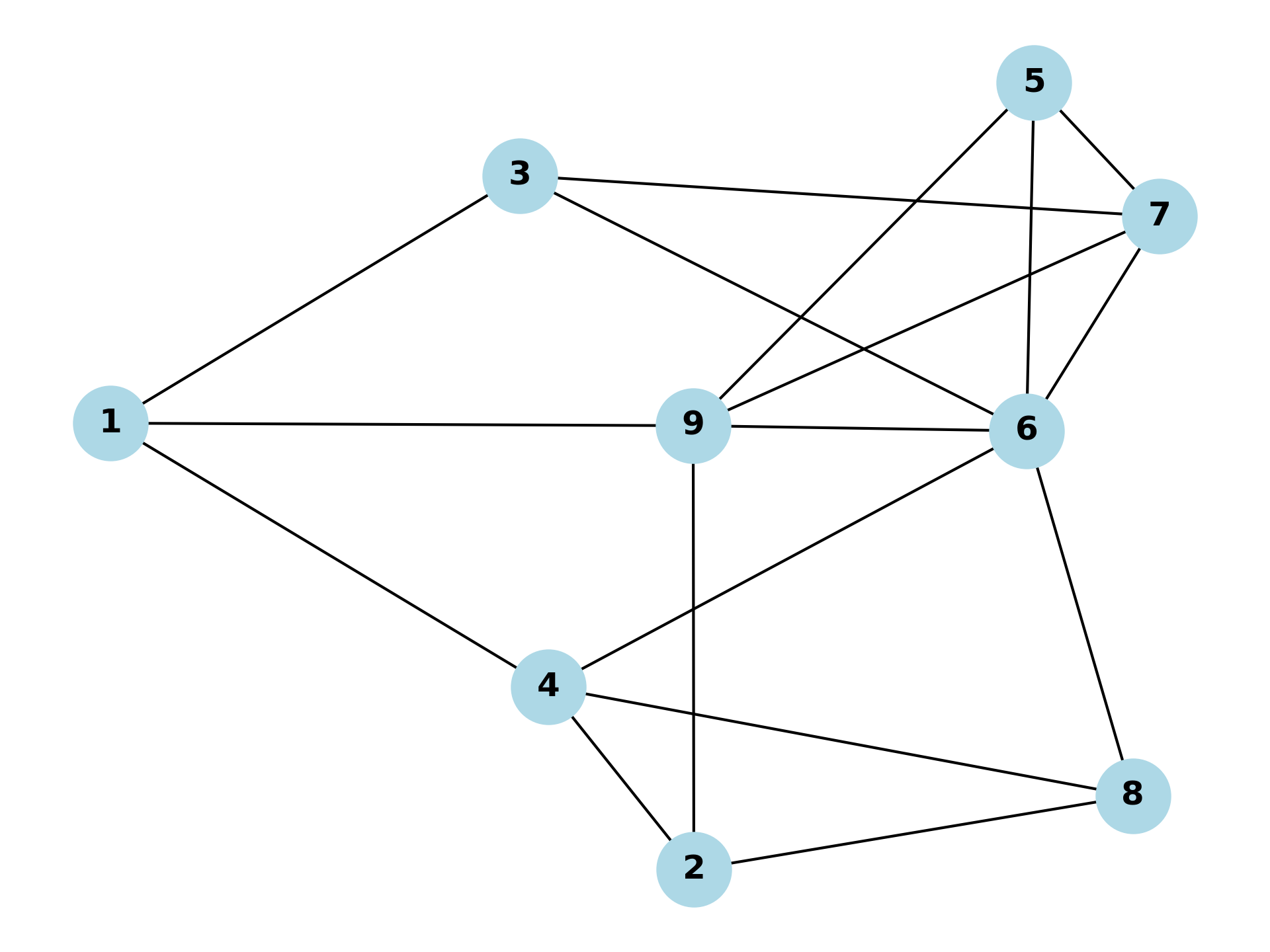}
		\caption{$n=9,\delta(G)=3,I'(G)=3$}
		\label{fig74}%文中引用该图片代��?
	\end{subfigure}
	\centering
	\begin{subfigure}{0.325\linewidth}
		\centering
		\includegraphics[width=1\linewidth]{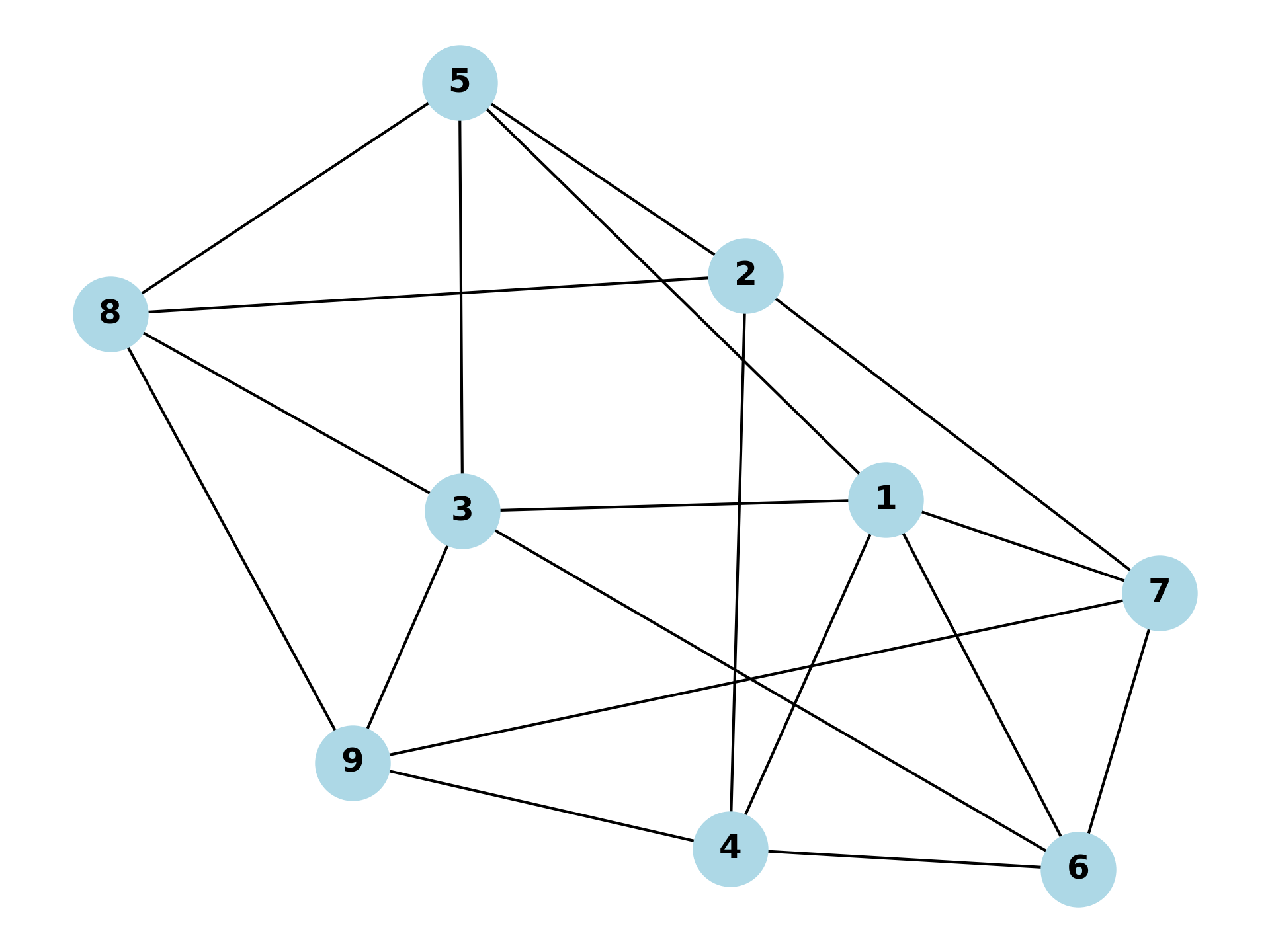}
		\caption{$n=9,\delta(G)=4,I'(G)=3$}
		\label{fig75}%文中引用该图片代��?
	\end{subfigure}
	\centering
	\begin{subfigure}{0.325\linewidth}
		\centering
		\includegraphics[width=1\linewidth]{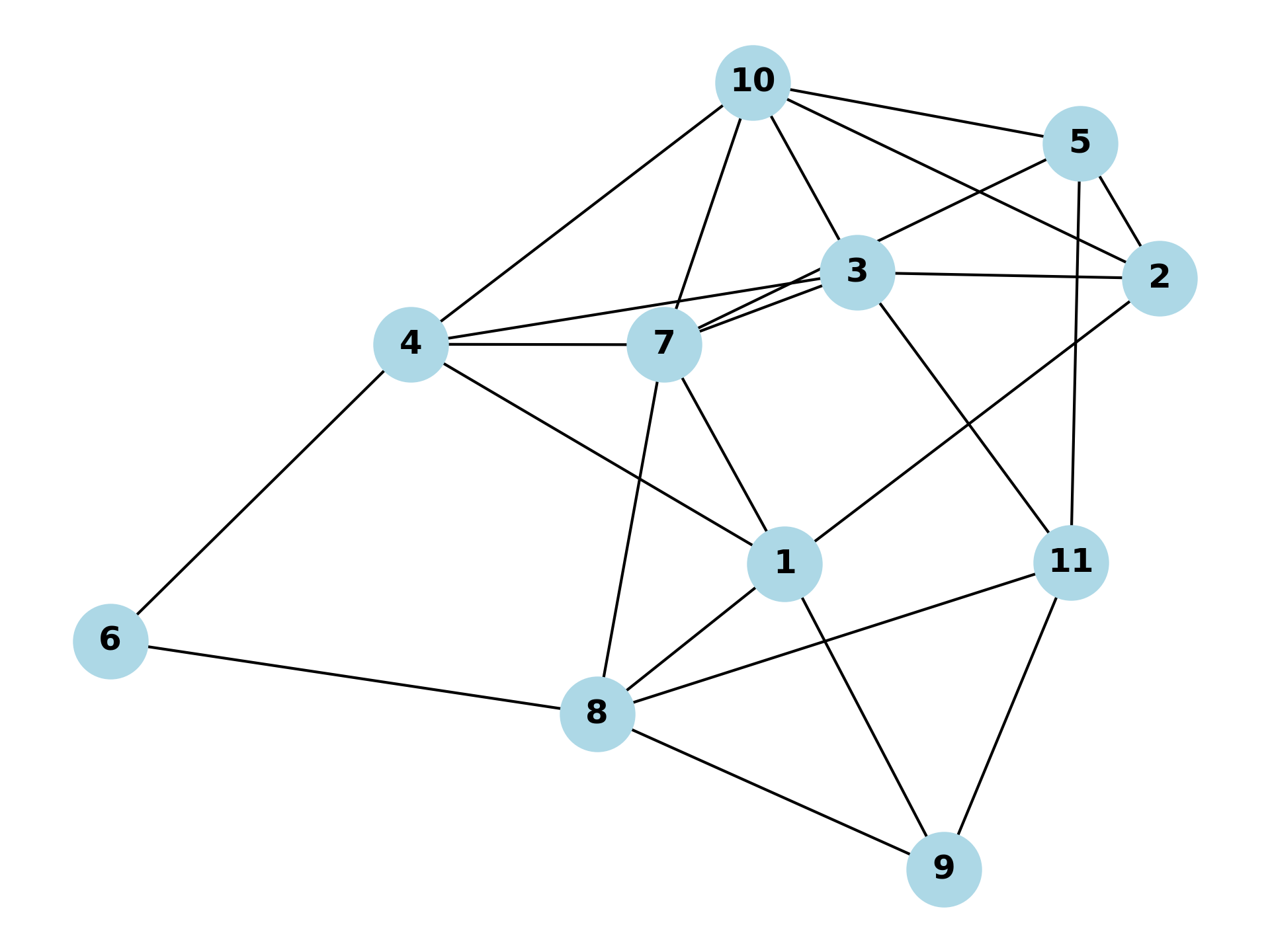}
		\caption{$n=11,\delta(G)=2,I'(G)=4$}
		\label{fig76}%文中引用该图片代��?
	\end{subfigure}

	\begin{subfigure}{0.325\linewidth}
		\centering
		\includegraphics[width=1\linewidth]{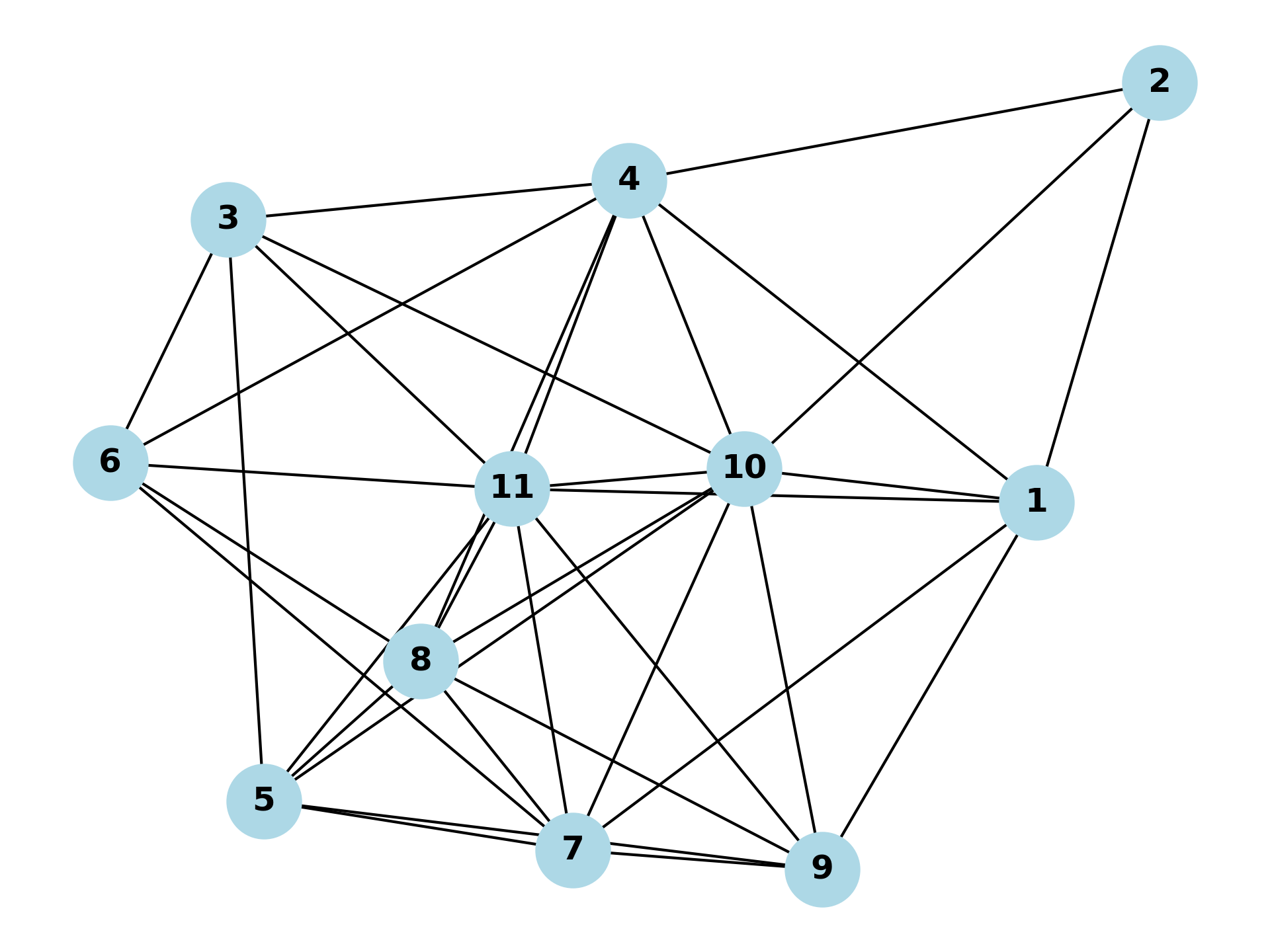}
		\caption{$n=11,\delta(G)=3,I'(G)=4$}
		\label{fig77}%文中引用该图片代��?
	\end{subfigure}
	\centering
	\begin{subfigure}{0.325\linewidth}
		\centering
		\includegraphics[width=1\linewidth]{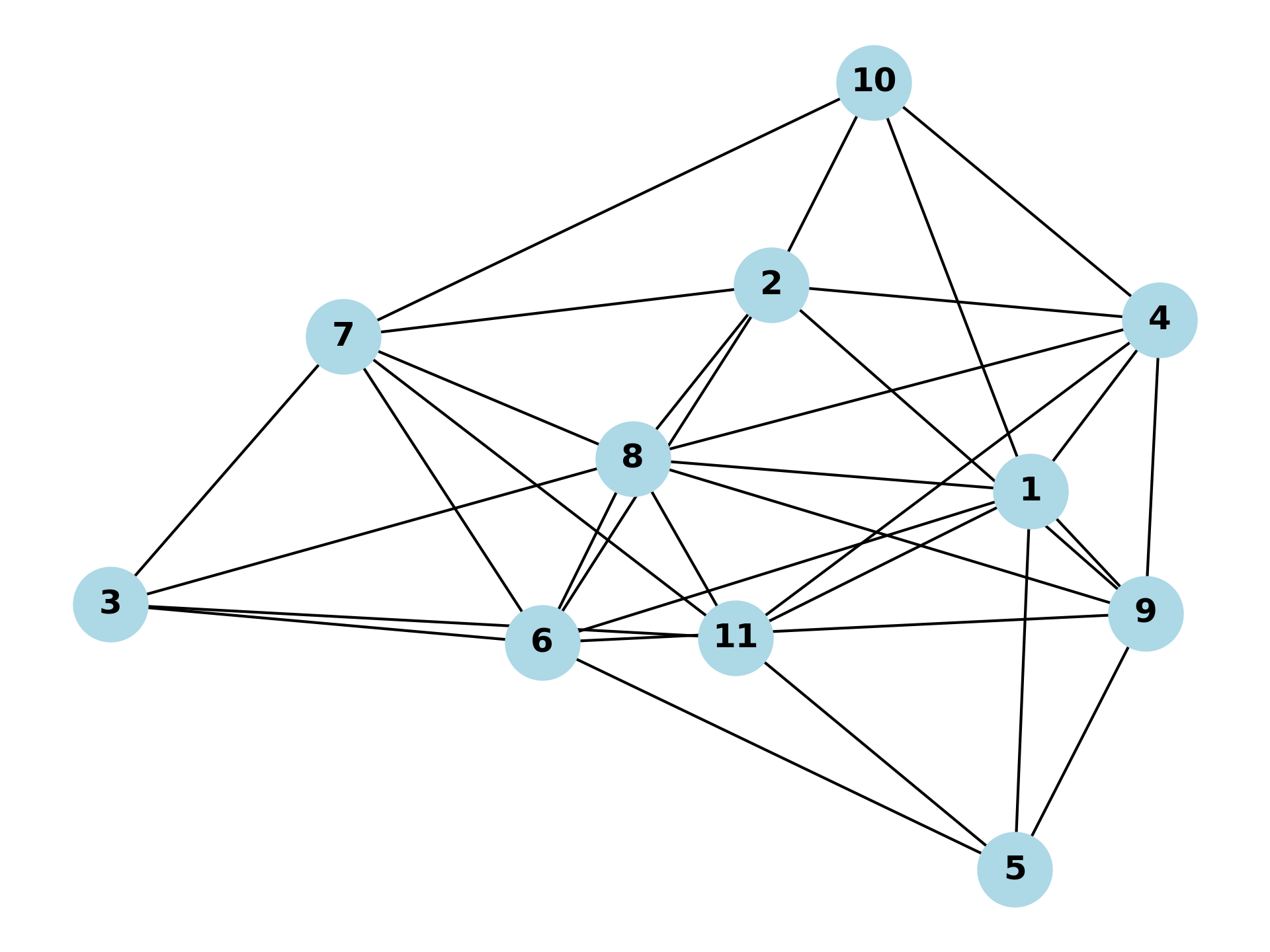}
		\caption{$n=11,\delta(G)=4,I'(G)=4$}
		\label{fig78}%文中引用该图片代��?
	\end{subfigure}
	\centering
	\begin{subfigure}{0.325\linewidth}
		\centering
		\includegraphics[width=1\linewidth]{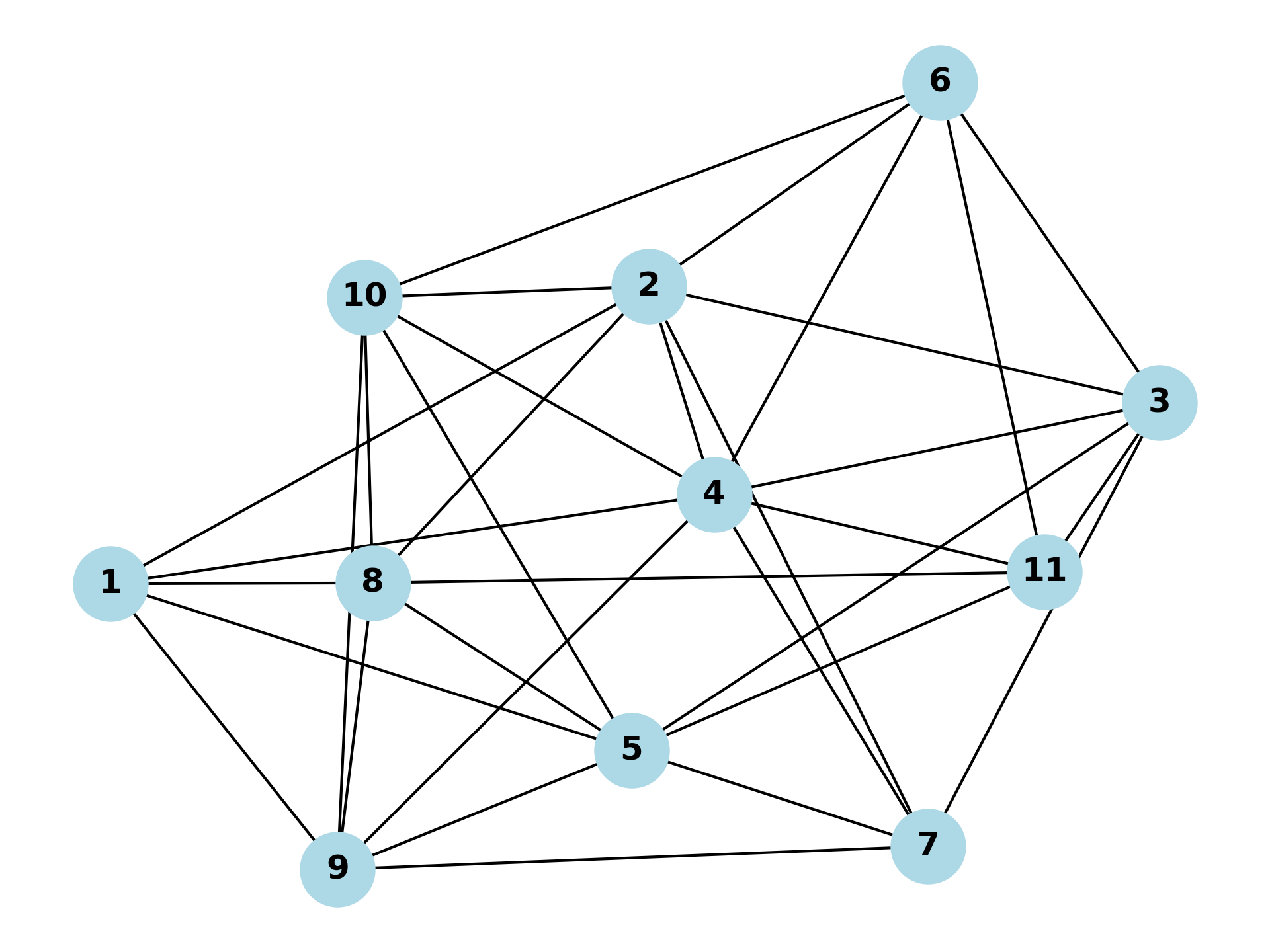}
		\caption{$n=11,\delta(G)=5,I'(G)=2.333$}
		\label{fig79}%文中引用该图片代��?
	\end{subfigure}
	\caption{Visualization of partial graph structures in Tab. \ref{tab3} and Tab. \ref{tab4}.}
	\label{fig7}
\end{figure}

\begin{table}[!htp]
    \centering
    \caption{The number of graph structures that satisfy the corresponding requirements (average and variance of thirty experiments) before the implementation of diversity enhancement in the algorithm.} \label{tab4}
    \begin{tabular}{ccllllll}
        \hline
        $n$  & $k$& \multicolumn{6}{c}{The number of graphs (average and variance)} \\
        \hline
			7		 &      2         &         36.7(3.742)&  33.5(3.667)  &  & & & \\
			9		 &      2         &         1.4(0.663) & 190.4(12.117) & 39.9(4.571) &  & &  \\
			11		 &      2         &         36.1(4721) & 110.8(11.956) & 90.9(5.682) & 3.6(1.356)    &  &  \\
			13		 &      2         &         3.6(2.244) & 214.5(17.67) & 210.9(10.290) & 53.5(8.321)  &3.0(1.612) &    \\
			15		 &      2         &         19.3(2.451) & 77.1(6.759)  &163.5(11.749) & 53.5(5.142)& 5.3(2.410)  &2.2(0.979)    \\
			15		 &      3         &         3.0(0.774) & 15.1(2.3) & 25.5(5.192);  &39.4(4.127)  &2.6(1.356)   &  \\
			15		 &      4         &        0.0(0.0)&  28.1(4.571) & 13.0(3.376) & 2.2(1.166)   &   &\\
        \hline
    \end{tabular}
 \end{table}
Fig. \ref{fig7} demonstrates the visualization of some graph structures that meet the requirements when $n$ is 7, 9, and 11, respectively, with $k=2$. It illustrates that our algorithm can provide multiple satisfactory graph structures in a relatively short period of time.

\subsection{Experimental analysis}
From the experimental setup in this paper, assuming that no identical individuals appear during the evolution process, the experiment can explore a maximum of 1000 graph structures, with each graph structure calculated for $I'(G)$ according to Alg. \ref{alg:pga}. Although this sacrifices the possibility of generating more graph structures that meet the conditions, it significantly speeds up the process. Compared to the traditional brute-force enumeration method, the time consumption is greatly reduced, as shown in Tab. \ref{tab:time}.

\begin{table}[!htp]
    \centering
    \caption{Compared to brute-force enumeration, the time consumption and results. All algorithms were run ten times, with the average values (standard deviation) recorded.} \label{tab:time}
    \begin{tabular}{cccccc}
        \hline
        $n$  & $k$&our time& $(\delta(G),I'(G))$&brute-force enumeration time&$(\delta(G),I'(G))$\\
        \hline
			6		 &      2         &  0.261(0.0051) &(2,4)  &7.740(0.021)&(2,4) \\
			7		 &      2         &  0.341(0.0027) &(2,5);(3,5)    &1134.4854(5.921)&(2,5);(3,5) \\

        \hline
    \end{tabular}
 \end{table}

Tab. \ref{tab:time} shows that our algorithm achieves the same results with significantly lower time costs. It is worth noting that we also experimented with brute-force enumeration for cases where the number of vertices is greater than or equal to 8, but due to its enormous time consumption, the results are not presented($O(2^{\frac{n(n-1)}{2}}2^n$)).

\section{Conclusion and discussion}
The design of the network structure is the most critical link in the entire network lifecycle, which is closely related to the subsequent maintenance, operation and management. Traditional network structure designing is routinely based on heuristic tricks, e.g., the network topology is scheduled in light of existing knowledge and expert experience. According to our existing knowledge, this paper is the first work to design the network structure thoroughly relying on graph-based parameters. Since fractional factors and isolated toughness variants are both celebrated indicators in the network, the parameter combination of $\delta(G)$ and $I'(G)$ obtained by the existence of fractional factors tremendously reflects the data transmission capacity and vulnerability of the network. We design an inverse model algorithm to infer the structure of the graph (network) $G$ with the help of $(\delta(G),I'(G))$ and the evolutionary computing strategy.

Although the algorithm given in this paper has been verified to be highly efficient in some specific scenarios, it still involves some defects. In the conventional scenario of network construction, the weighted information is endowed, i.e., we know a weight function $\xi: V\times V\to \Bbb R$ which maps each pair of vertices to a real number (for instance, $\xi$ can be a distance function, connectivity, or kernel function between two vertices). Then the scale of the weight determines whether to construct an edge for each pair of vertices. However, the setting in this article ignores the weights of vertex pair, because the isolated toughness variant is closely related to the vertex cut set, which is a topology-based parameter and extraneous with weight. This is the crux reason why the isolated toughness variant has no corresponding definition in the weighted graph setting. In this case, the implicit meaning lies in that the relationship between any pair of vertices is evenly treated. Therefore, the network structure constructed by the approach in this article may have a significant gap with the ideal network. How to incorporate the relationship between sites (such as physical distance) as prior knowledge into the inverse model algorithm needs to be further studied.

Finally, we discuss the following problems from theoretical perspectives. For a given graph $G$, let
$$S_{G}^{\dagger}=\mathop{\arg\min}_{S\subseteq V(G),i(G-S)\ge2}\frac{|S|}{i(G-S)}$$
and
$$S_{G}^{\ddagger}=\mathop{\arg\min}_{S\subseteq V(G),i(G-S)\ge2}\frac{|S|}{i(G-S)-1}.$$
It implies $S_{G}^{\dagger}$ and $S_{G}^{\ddagger}$ are the vertex subsets to achieve the isolated toughness and isolated toughness variant, respectively. From the discussion given in this paper we know that there is an inherent difference between $I(G)$ and $I'(G)$. However, so far, all graph families meet $S_{G}^{\dagger}=S_{G}^{\ddagger}$. For instance, for $G_{l}=K_{l-1}\vee (lK_{k})$, we have $S_{G_{l}}^{\dagger}=S_{G_{l}}^{\ddagger}$ for any  $l\in\Bbb N-\{1\}$. It is natural to conduct the following problem.
\begin{problem}\label{problem2} Whether $|S_{G}^{\dagger}|=|S_{G}^{\ddagger}|$ for any graph $G$? If the answer is negative, then can we provide a graph family $\mathcal{G}$ satisfying $|S_{G}^{\dagger}|\ne|S_{G}^{\ddagger}|$ if $G\in\mathcal{G}$.
\end{problem}
Till now, we have only the following clue for Problem \ref{problem2}.
\begin{proposition}\label{proposition1} Let $G$ be a graph with $|S_{G}^{\dagger}|\ne |S_{G}^{\ddagger}|$. Then $|S_{G}^{\ddagger}|>|S_{G}^{\dagger}|$ and $i(G-S_{G}^{\ddagger})>i(G-S_{G}^{\dagger})$.
\end{proposition}
{\bf Proof of Proposition \ref{proposition1}.} According to the definition of $S_{G}^{\dagger}$ and $S_{G}^{\ddagger}$, we deduce $\frac{|S_{G}^{\dagger}|}{i(G-S_{G}^{\dagger})}\le\frac{|S_{G}^{\ddagger}|}{i(G-S_{G}^{\ddagger})}$ and $\frac{|S_{G}^{\ddagger}|}{i(G-S_{G}^{\ddagger})-1}\le\frac{|S_{G}^{\dagger}|}{i(G-S_{G}^{\dagger})-1}$. Hence
$$|S_{G}^{\ddagger}|(i(G-S_{G}^{\dagger})-1)\le|S_{G}^{\dagger}|(i(G-S_{G}^{\ddagger})-1),$$
which reveals
$$\frac{|S_{G}^{\ddagger}|}{i(G-S_{G}^{\ddagger})}+\frac{|S_{G}^{\dagger}|-|S_{G}^{\ddagger}|}{i(G-S_{G}^{\ddagger})i(G-S_{G}^{\dagger})}\le\frac{|S_{G}^{\dagger}|}{i(G-S_{G}^{\dagger})}.$$
Hence, the term ``$\frac{|S_{G}^{\dagger}|-|S_{G}^{\ddagger}|}{i(G-S_{G}^{\ddagger})i(G-S_{G}^{\dagger})}$'' must be negative, which implies the desired proposition. \hfill $\Box$

Due to the fact that computational toughness is an NP-hard problem (see Bauer et al. \cite{Bauer90}), intuitively speaking,  the computation of $I(G)$ or $I'(G)$ is also an NP-hard problem. However, strictly speaking, we have not found any literature discussing the computational complexity of isolated toughness and its variant. That is to say, this seemingly obvious problem has not been rigorously proven yet. It naturally leads us to  the following problem.
\begin{problem}\label{problem2} The computing of isolated toughness {\rm(}variant{\rm)} is an NP-hard problem?
\end{problem}

\section*{Availability of data and materials}
The code in this paper is made public on https://github.com/AizhEngHN/Inverse-model-for-network-construction-G-I-G-G.

\section*{Competing interests}
We declare that none of the authors have any competing interests in the paper.

\section*{Acknowledgment}

We thank the reviewers for their constructive comments in
improving the quality of this paper. This work has been partially
supported by National Science Foundation of China (Nos. 12161094, 11871270 and 11931006).

\end{document}